

  \newcount\fontset
  \fontset=2
  \def\dualfont#1#2#3{\font#1=\ifnum\fontset=1 #2\else#3\fi}
  \dualfont\eightrm {cmr8} {cmr7}
  \dualfont\eightmit {cmmi8} {cmmi7}
  \dualfont\eightit {cmti8} {cmti10}
  \dualfont\eightbf {cmbx8} {cmr7}
  \dualfont\tensc {cmcsc10} {cmr10}
  \dualfont\twelvess {cmss10 scaled \magstep2} {cmr10 scaled
\magstep2}

  \magnification=\magstep1
  \nopagenumbers
  \voffset=2\baselineskip
  \advance\vsize by -\voffset
  \headline{\ifnum\pageno=1 \hfil\else\ifodd\pageno
    \tensc\hfil bunce-deddens algebras and partial automorphisms
    \hfil\folio
    \else\tensc\folio\hfil ruy exel\hfil\fi\fi}

  
  \def\<{\langle}
  \def\>{\rangle}
  \def\:{\colon}
  \def\+{\oplus}
  \def\*{\otimes}
  \def\arw{\rightarrow}
  \def\cstar{$C^*$}
  \def\eightcstar{{\eightmit C\kern1pt}\raise.4ex\hbox{\fiverm *}}
  \def\C{{\bf C}}
  
  \def\N{{\bf N}}
  
  \def\th{\null^{\hbox{\sevenrm th}}}

  \def\sameauthor{\underbar{\hbox to 1.5truecm{\hfil}}}
  \def\square{\hbox{$\sqcap\!\!\!\!\sqcup$}}

  \newcount\bibno \bibno=0
  \def\newbib#1{\advance\bibno by 1 \edef#1{\number\bibno}}
  \def\stdbib#1#2#3#4#5#6#7#8{\smallskip \item{[#1]} #2, ``#3'',
    {\sl#4} {\bf#5} (#6), #7--#8.}
  \def\bib#1#2#3#4{\smallskip \item{[#1]} #2, ``#3'', {#4.}}
  \def\cite#1{{\rm[\bf #1\rm]}}
  \def\speccite#1#2{{\rm[\bf #1\rm, #2]}}
  \outer\def\beginsection#1\par{\vskip0pt plus.1\vsize\penalty-250
    \vskip0pt plus-.1\vsize\bigskip\vskip\parskip
    \message{#1}\centerline{\tensc#1}\nobreak\bigskip}
  \def\statement#1#2#3{\medbreak\noindent{\bf#1.\enspace
    #2.\enspace}{\sl #3}\medbreak}
  
  \def\endproof{\ifmmode\eqno\square\else\hfill\square\medbreak\fi}
  \def\today{\ifcase\month\or January\or February\or March\or April\or
    May\or June\or July\or August\or September\or October\or
    November\or December\fi\ \number\day, \number\year}

  \def\ltwo{l_2}
  \def\linf{l_\infty}
  \def\BDT{A}
  \def\BD{B}
  \def\act{\alpha}
  \def\Ca{{\bf K}}
  \def\O{{\cal O}}
  \def\K{{\cal K}}
  \def\B{{\cal B}}

  \newbib{\Archbold}
  \newbib{\Blackadar}
  \newbib{\BK}
  \newbib{\BDI}
  \newbib{\BDII}
  \newbib{\ER}
  \newbib{\NP}
  \newbib{\AF}
  \newbib{\Green}
  \newbib{\Ped}
  \newbib{\Power}
  \newbib{\Putnam}
  \newbib{\Riedel}

  \null
  \begingroup
  \def\r{\rightline}
  \def\c{\centerline}

  \vskip -\voffset \vskip -2truecm \fiverm \baselineskip=3ex

                                                  \r{UNM--RE--005}
                                              \r{February 8, 1993}
                                                \r{Printed \today}

  \vskip 1.5truecm \twelvess \baselineskip=3ex

	  \c{The Bunce-Deddens Algebras as Crossed Products}
		     \c{by Partial Automorphisms}

 \bigskip \tensc \baselineskip=3ex

	       \c{Ruy Exel\footnote{\raise.5ex\hbox{*}}
       {\eightrm On leave from the University of S\~ao Paulo.}}

  \bigskip \eightit \baselineskip=3ex

	     \c{Department of Mathematics and Statistics,
		      University of New Mexico}
		  \c{Albuquerque, New Mexico 87131}

  \bigskip\bigskip\eightrm\baselineskip=3ex
  \midinsert \narrower We describe both the Bunce-Deddens
\eightcstar-algebras and their Toeplitz versions, as crossed products
of commutative \eightcstar-algebras by partial automorphisms.  In the
latter case, the commutative algebra has, as its spectrum, the union
of the Cantor set and a copy the set of natural numbers {\eightbf N},
fitted together in such a way that {\eightbf N} is an open dense
subset. The partial automorphism is induced by a map that acts like
the odometer map on the Cantor set while being the translation by one
on {\eightbf N}.  From this we deduce, by taking quotients, that the
Bunce-Deddens \eightcstar-algebras are isomorphic to the (classical)
crossed product of the algebra of continuous functions on the Cantor
set by the odometer map. \endinsert \endgroup

  \beginsection 1. Introduction

Recall from \cite{\BDI} that a weighted shift operator is a bounded
operator on $\ltwo = \ltwo(\N)$ given, on the canonical basis
$\{e_n\}_{n=0}^\infty$, by $S_a(e_n) = a_{n+1}e_{n+1}$, where the
weight sequence $a=\{a_n\}_{n=1}^\infty$, is a bounded sequence of
complex numbers. A weighted shift is said to be $p$-periodic if its
weights satisfy $a_n = a_{n+p}$ for all $n$.

Given a strictly increasing sequence $\{n_k\}_{k=0}^\infty$ of
positive integers, such that $n_k$ divides $n_{k+1}$ for all $k$, the
Bunce-Deddens-Toeplitz \cstar-algebra $\BDT = \BDT(\{n_k\})$ is
defined to be the \cstar-algebra of operators on $\ltwo$, generated by
the set of all $n_k$-periodic weighted shifts, for all $k$. These
algebras were first studied by Bunce and Deddens in \cite{\BDII}. It
was observed by them that the algebra $\K$, of compact operators on
$\ltwo$, is contained in $\BDT$ and that the quotient $\BDT/\K$ is a
simple \cstar-algebra. The latter became known as the Bunce-Deddens
\cstar-algebra and has been extensively studied (see, for example,
  \cite{\Archbold},
  \speccite{\Blackadar}{10.11.4},
  \cite{\BK},
  \cite{\ER},
  \speccite{\Green}{p.~248} ,
  \cite{\Power},
  \cite{\Putnam},
  \cite{\Riedel}).
  We shall denote these algebras by $\BD(\{n_k\})$ or simply by $\BD$,
if the weight sequence is understood.

The goal of the present work is to describe both $\BDT(\{n_k\})$ and
$\BD(\{n_k\})$ as the crossed product of commutative \cstar-algebras
by partial automorphisms \cite{\NP}, in much the same way as we have
described general AF-algebras \cite{\AF} as partial crossed products.

In the case of $\BDT(\{n_k\})$, we shall see that it is given by a
curious (partial) dynamical system consisting of a topological space
$X$ which can be thought of as a compactification of the (discrete)
space $\N$ of natural numbers, the complement of $\N$ in $X$ being
homeomorphic to the Cantor set $\Ca$. The transformation $f$ of $X$,
by which the partial crossed product is taken, leaves both $\N$ and
$\Ca$ invariant. Its behavior on $\N$ is that of the translation by
one, while the action on $\Ca$ is by means of the odometer map (see,
for example, \cite{\Putnam}) which is defined as follows. Given a
sequence $\{q_k\}_{k=0}^\infty$ of positive integers (below we shall
use $q_k = n_{k+1}/n_k$), consider the Cantor set, as given by the
model
  $$\Ca = \prod_{j=0}^{\infty}\{0,1,\ldots,q_j-1\}.$$
  The odometer map is the map\quad $\O\: \Ca \arw \Ca$, \quad given by
formal addition of $(1,0, \ldots)$ with carry over to the right.  Note
that \quad $\O(q_0-1, q_1-1, \ldots) = (0, 0, \ldots )$ \quad since
the carry over process, in this case, extends all the way to infinite.
For further reference let us call by the name of ``partial odometer''
the restriction of $\O$ to a map from $X-\{(q_0-1, q_1-1, \ldots)\}$
to $X-\{(0, 0, \ldots)\}$ so that, for this map, the carry over
process always terminates in finite time.

As already mentioned, $\BD(\{n_k\})$ is the quotient of
$\BDT(\{n_k\})$ by $\K$.  But $\K$ can be seen to correspond to the
restriction of the above dynamical system to $\N$ (see \cite{\NP}).
So, we deduce that the Bunce-Deddens algebras $\BD(\{n_k\})$ are
isomorphic to the crossed product of the Cantor set by the odometer
map. This result is already well known
  \speccite{\Blackadar}{10.11.4}
  but it is interesting to remark how little bookkeeping is necessary
to deduce it from the machinery of partial automorphisms \cite{\NP}.
  Moreover, this should be compared with \cite{\AF}, Theorem 3.2,
according to which UHF-algebras correspond to the crossed product of
the Cantor set by the partial odometer map. One therefore obtains a
crystal clear picture of the fact, already noticed by Bunce and
Deddens, that UHF-algebras sit as subalgebras of $\BD(\{n_k\})$ (see
also \cite{\Putnam}).

  \beginsection 2. Circle Actions

For each $z$ in the unit circle $S^1 = \{w\in \C \: |w| = 1\}$ let,
$U_z$ denote the diagonal unitary operator on $\ltwo$, given by
$U_z(e_n) = z^ne_n$. If $S$ is any weighted shift, it is easy to see
that $U_zSU_z^{-1} = zS$. Thus, denoting by $\act_z$ the inner
automorphism of $\B(\ltwo)$ given by conjugation by $U_z$, one finds
that $\BDT$ is invariant under $\act_z$. Moreover, one can see that
this gives a continuous action of $S^1$ on $\BDT$, in the sense of
\cite{\Ped}, 7.4.1 (even though the corresponding action is not
continuous on $\B(\ltwo)$).

Let us denote the fixed point subalgebra by $\BDT_0$. It is easy to
see that $\BDT_0$ consists precisely of the operators in $\BDT$ which
are diagonal with respect to the basis $\{e_n\}$. Now, given that the
\cstar-algebra of (bounded) diagonal operators is isomorphic to $\linf
= \linf(\N)$, it is convenient to view $\BDT_0$ as a subalgebra of
$\linf$. For the purpose of describing $\BDT_0$, observe that, since
$\K$ is contained in $\BDT$, it follows that $c_0$ (the subalgebra of
$\linf$ formed by sequences tending to zero) is contained in $\BDT_0$.
Carrying this analysis a bit further one can prove that $\BDT_0 = c_0
\oplus D$ where $D$ is the subalgebra of $\linf$ generated by all
$n_k$-periodic sequences, for all $k$.

This decomposition is useful in determining the spectrum of $\BDT_0$.
Note, initially, that spectrum of $c_0$ is homeomorphic to $\N$ (with
the discrete topology) while the spectrum of $D$ is the Cantor set,
here denoted $\Ca$. This said, one has that the spectrum of $\BDT_0$
can be described, at least in set theoretical terms, as the union $X =
\N \cup \Ca$.  Moreover, since $c_0$ is an essential ideal in
$\BDT_0$, one sees that $\N$ is an open dense subset of $X$. To better
grasp the entire topology of $X$ we need a more precise notation.
Assume, without loss of generality, that $n_0=1$ and let $q_k =
n_{k+1}/n_k$ for $k\geq 0$.
  Any integer $n$ with $0 \leq n < n_k$ has a unique representation as
  $$ n= \sum_{j=0}^{k-1}\beta^{(n)}_jn_j$$
  where $0 \leq \beta^{(n)}_j < q_j$. Here the $\beta^{(n)}_j$ play
the role of digits in a decimal-like representation, except that the
base varies along with the position of each digit. Accordingly, we let
$\beta^{(n)} = (\beta^{(n)}_0, \ldots, \beta^{(n)}_{k-1})$ be the
corresponding notation for $n$ (which we shall use interchangeably
without further warning).  When convenient, we shall also view
$\beta^{(n)}$ as an element of the set
  $$\Ca_k = \prod_{j=0}^{k-1} \{ 0, 1, \ldots, q_j-1\}.$$
  For each $k$ and each $\beta$ in $\Ca_k$, we denote by $e_\beta$ the
$n_k$-periodic sequence (thus an element of $D$) given by
  $$e_\beta(n) =
  \left\{\matrix{1 & \hbox{if\quad} n\equiv \beta\pmod{n_k} \cr
                 0 & \hbox{otherwise}\hfill}
  \right.$$ Note that the length of $\beta$ determines which $n_k$
should be used in the above definition.
  Clearly the set $\{e_\beta \: \beta\in \bigcup_{k=1}^\infty\Ca_k\}$
generates $D$. Making use of the notation introduced, we shall adopt
for the Cantor set, the model
  $$\Ca = \prod_{j=0}^{\infty} \{ 0, 1, \ldots, q_j-1\}$$
  so that, once we view $D$ as the algebra of continuous functions on
$\Ca$ via the Gelfand transform, the support of the $e_\beta$ form a
basis for the topology of $\Ca$. In fact, the support of $e_\beta$ is
precisely the set of elements $\gamma = (\gamma_i)$ in $\Ca$ such that
$\gamma_j = \beta_j$ for all $j=0, \ldots, k-1$ (assuming that $\beta$
is in $\Ca_k$). That is, $\gamma$ is in the support of $\beta$ if and
only if its initial segment coincides with $\beta$.

Considering now, the whole of $\BDT_0$, note that it is generated by
the set of all idempotents $p$ which, viewed as elements of $\linf$,
have one of the two following forms: either it has a finite number of
non-zero coordinates (in which case $p$ is in $c_0$), or it coincide
with some $e_\beta$, except for finitely many coordinates. The set of
such idempotents is closed under multiplication, which therefore
implies that their support, in the spectrum $X$ of $\BDT_0$, form a
basis for the topology of $X$.  With this we have precisely described
the topology of $X$:

\statement{2.1}{Theorem}{The spectrum of $\BDT_0$ consists of the
union of the Cantor set\/ $\Ca$ and a copy of the set of natural
numbers\/ $\N$.  Each element of\/ $\N$ is an isolated point and a
fundamental system of neighborhoods of a point $\gamma = (\gamma_i)$
in\/ $\Ca$ consists of the sets $V_{k}$ defined to be the union of the
sets
  $$\{ \zeta \in \Ca\: \zeta_i = \gamma_i, i < k\}$$
  and
  $$\{ n \in \N\: n \geq k {\quad\rm and\quad} \beta^{(n)}_i =
\gamma_i, i < k\}.$$}
  Note the interesting interplay between the digital representation of
the natural numbers on one hand, and of elements of the Cantor set, on
the other.

  \beginsection 3. The Main Result

Recall from \cite{\NP}, Theorem 4.21, that a regular semi-saturated
action of $S^1$ on a \cstar-algebra, causes it to be isomorphic to the
covariance algebra of a certain partial automorphism of the fixed
point subalgebra.
  In the case of the action $\act$ of $S^1$ on $\BDT$, described
above, it is very easy to see that it is regular and semi-saturated.
The semi-saturation follows immediately, since every weighted shift
belongs to the first spectral subspace of $\act$, henceforth denoted
$\BDT_1$.
  The fact that $\act$ is regular depends on the existence of an
isomorphism
  \quad $\theta\: \BDT_1^* \BDT_1 \rightarrow \BDT_1 \BDT_1^*$\quad
  and a linear isometry
  \quad $\lambda\: \BDT^*_1 \rightarrow \BDT_1 \BDT_1^*$\quad
  such that, for $x, y \in \BDT_1$, $a \in \BDT_1^* \BDT_1$ and $b \in
\BDT_1 \BDT_1^*$
  \medskip
  \itemitem{(i)} $\lambda (x^* b) = \lambda(x^*) b$
  \medskip
  \itemitem{(ii)} $\lambda (ax^*) = \theta (a) \lambda (x^*)$
  \medskip
  \itemitem{(iii)} $\lambda (x^*)^* \lambda (y^*) = xy^*$
  \medskip
  \itemitem{(iv)} $\lambda (x^*) \lambda (y^{*})^{*} = \theta (x^*
y).$

\medskip\noindent See \cite{\NP}, 4.3 and 4.4 for more information.
It is easy to see that $\theta$ and $\lambda$, given by
  $\theta(a) = S a S^*$ and
  $\lambda(x^*) = Sx^*$,
  satisfy the desired properties, where $S$ denotes the unilateral
(unweighted) shift.

Note that, in the present case, $\BDT_1^*\BDT_1 = \BDT_0$ while
$\BDT_1\BDT_1^*$ is the ideal of $\BDT_0$ formed by all sequences for
which the first coordinate vanishes. That is, under the standard
correspondence between ideals and open subsets of the spectrum,
$\BDT_1\BDT_1^*$ corresponds to $X-\{0\}$. The isomorphism
  \quad $\theta\: C(X) \rightarrow C(X-\{0\})$\quad
  therefore induces a homeomorphism
  \quad $f\: X \rightarrow X-\{0\}$\quad
  which we would now like to describe.

If an integer $n$ is thought of as an element of $X$, as seen above,
then the element $\delta_n$ of $\linf$,
  represented by the sequence having the $n\th$ coordinate equal to
one and zeros everywhere else,
  corresponds to the characteristic function of the singleton $\{n\}$
and, given that $\theta(\delta_n) = \delta_{n+1}$, we see that $f(n) =
n+1$.

  Now, if $\gamma=(\gamma_i)$ is in $\Ca$, let $\gamma|_k$ be the
$k\th$ truncation of $\gamma$, that is $\gamma|_k = (\gamma_0, \ldots,
\gamma_{k-1})$ so that we can consider $e_{\gamma|_k}$, as defined
above.  Also let $f_k$ be the element of $c_0$ given by $f_k = (1,
\ldots, 1, 0, \ldots)$
  where the last ``1'' occurs in the position $k-1$, counting from
zero.
  The support of the Gelfand transform of the idempotent element
$(1-f_k)e_{\gamma|_k}$ is precisely the set $V_k$ referred to in 2.1.
So, as $k$ varies, the intersection of these sets is precisely the
singleton $\{\gamma\}$.  Therefore, to find out what $f(\gamma)$
should be, it is enough to look for the intersection of the supports
of the Gelfand transforms of the elements
  $$\theta\bigl((1-f_k)e_{\gamma|_k}\bigr) = S(1-f_k) e_{\gamma|_k}
S^*.$$

The reader is now invited to verify, using this method, that the
effect that $f$ has on $\gamma$ is precisely the effect of the
odometer map. One should exercise special attention to check that the
above method does indeed give
  \quad $f(q_0-1, q_1-1, \ldots) = (0, 0, \ldots)$,\quad
  in contrast with the partial automorphisms that produce UHF-algebras
\cite{\AF}, since $(q_0-1, q_1-1, \ldots)$ is removed from the domain
of the maps considered there.
  Summarizing our findings so far, we have:

\statement{3.1}{Theorem}{The Bunce-Deddens-Toeplitz \cstar-algebra
$\BDT(\{n_k\})$ is isomorphic to the crossed product of\quad $C(X)$ by
the partial automorphism
  \quad $\theta\: C(X) \rightarrow C(X-\{0\})$\quad
  induced by the (inverse of the) map
  \quad $f\: X \rightarrow X-\{0\}$\quad
  acting like the odometer on\/ $\Ca$ and like translation by one on\/
$\N$.}

Recall that the Bunce-Deddens \cstar-algebras were defined to be the
quotient $\BD = \BDT/\K$. If one identifies $\BDT$ with the crossed
product above, it is easy to see that the ideal $\K$ corresponds to
the crossed product of $c_0$ by the corresponding restriction of
$\theta$. Therefore, the quotient can be described as the (classical)
crossed product of the Cantor set by the odometer map. That is:

\statement{3.2}{Theorem}{The Bunce-Deddens \cstar-algebra
$\BD(\{n_k\})$ is isomorphic to the crossed product of $C(\Ca)$ by the
automorphism induced by the odometer map.}

  \bigbreak
  \centerline{\tensc References}
  \nobreak\medskip
  \frenchspacing

  \stdbib{\Archbold}
  {R. J. Archbold}
  {An averaging process for \cstar-algebras related to weighted
shifts}
  {Proc. London Math. Soc.} {35} {1977} {541} {554}

  \bib{\Blackadar}
  {B. Blackadar}
  {$K$-theory for operator algebras}
  {MSRI Publications, Springer--Verlag, 1986}

  \stdbib{\BK}
  {B. Blackadar and A. Kumjian}
  {Skew products of relations and the structure of simple
\cstar-algebras}
  {Math. Z.} {189} {1985} {55} {63}

  \stdbib{\BDI}
  {J. W. Bunce and J. A. Deddens}
  {\cstar-algebras generated by weighted shifts}
  {Indiana Univ. Math. J.} {23} {1973} {257} {271}

  \stdbib{\BDII}
  {\sameauthor}
  {A family of simple \cstar-algebras related to weighted shift
operators}
  {J. Funct. Analysis} {19} {1975} {13} {24}

  \stdbib{\ER}
  {E. E. Effros and J. Rosenberg}
  {\cstar-algebras with approximate inner flip}
  {Pacific J. Math.} {77} {1978} {417} {443}

  \bib{\NP}
  {R. Exel}
  {Circle actions on \cstar-algebras, partial automorphisms and a
  generalized Pimsner--Voiculescu exact sequence}
  {preprint, University of New Mexico, 1992}

  \bib{\AF}
  {\sameauthor}
  {Approximately finite \cstar-algebras and partial automorphisms}
  {preprint, University of New Mexico, 1992}

  \stdbib{\Green}
  {P. Green}
  {The local structure of twisted covariance algebras}
  {Acta Math.} {140} {1978} {191} {250}

  \bib{\Ped}
  {G. K. Pedersen}
  {\cstar-Algebras and their Automorphism Groups}
  {Academic Press, 1979}

  \stdbib{\Power}
  {S. C. Power}
  {Non-self-adjoint operator algebras and inverse systems of
simplicial complexes}
  {J. reine angew. Math.} {421} {1991} {43} {61}

  \stdbib{\Putnam}
  {I. F. Putnam}
  {The \cstar-algebras associated with minimal homeomorphisms of the
  Cantor set}
  {Pacific J. Math.}
  {136} {1989} {329} {353}

  \stdbib{\Riedel}
  {N. Riedel}
  {Classification of the \cstar-algebras associated with minimal
rotations}
  {Pacific J. Math.} {101} {1982} {153} {161}

  \bye